\documentclass{agtart_a}
\pdfoutput=1
\usepackage{graphicx}

\title{On the Kontsevich integral of Brunnian links}

\author{Kazuo Habiro}
\givenname{Kazuo}
\surname{Habiro}
\address{Research Institute for Mathematical Sciences\\
 Kyoto University\\\newline
 Kyoto 606-8502\\Japan}
\email{habiro@kurims.kyoto-u.ac.jp}
\urladdr{}

\author{Jean-Baptiste Meilhan}
\givenname{Jean-Baptiste}
\surname{Meilhan}
\email{meilhan@kurims.kyoto-u.ac.jp}
\urladdr{}

\volumenumber{6}
\issuenumber{}
\publicationyear{2006}
\papernumber{50}
\startpage{1399}
\endpage{1412}

\doi{}
\MR{}
\Zbl{}


\received{10 January 2006}
\revised{3 July 2006}
\accepted{7 July 2006}
\published{25 September 2006}
\publishedonline{25 September 2006}
\proposed{}
\seconded{}
\corresponding{}
\editor{}
\version{}

\arxivreference{math.GT/0605312}

\AtBeginDocument{\let\bar\wbar\let\tilde\wtilde\let\hat\what}
\def\SetFigFont#1#2#3#4#5{\small}%
\def\adjustlabel<#1,#2>#3{\smash{\rlap{\kern #1 \raise #2\hbox{#3}}}}
\newcommand\zzzcolon {\co}

\makeatletter
\def\cnewtheorem#1[#2]#3{\newtheorem{#1}{#3}[section]
\expandafter\let\csname c@#1\endcsname\c@theorem}
\makeatother

\theoremstyle{plain}
\newtheorem{theorem}{Theorem}[section]
\cnewtheorem{lemma}[theorem]{Lemma}
\cnewtheorem{proposition}[theorem]{Proposition}
\cnewtheorem{corollary}[theorem]{Corollary}
\theoremstyle{definition}
\cnewtheorem{definition}[theorem]{Definition}
\theoremstyle{remark}
\cnewtheorem{remark}[theorem]{Remark}
\newtheorem*{acknowledgments}{Acknowledgments}
\numberwithin{equation}{section}

\newcommand\Br{\operatorname{Br}}
\newcommand\BrJ{\Br(\bJ_{2n}(n+1))}
\newcommand\modZ {\mathbb{Z}}

\newcommand\modQ {\mathbb{Q}}
\newcommand\modL {\mathcal{L}}

\newcommand\bJ{\bar{J}}

\newcommand\modA {{\mathcal A}}

\newcommand\modB {{\mathcal B}}
\newcommand\modP {{\mathcal P}}

\newcommand\modone {{\mathbf 1}}
\newcommand\xto[1]{\overset{#1}{\longrightarrow}}

\newcommand\simeqto{\xto{\simeq}}

\newcommand\sumss{\sum_{\sigma ,\sigma '\in S_{n-1}}}

\begin{document}

\begin{asciiabstract}
The purpose of the paper is twofold.  First, we give a short proof
  using the Kontsevich integral for the fact that the restriction of
  an invariant of degree 2n to (n+1)-component Brunnian links can be
  expressed as a quadratic form on the Milnor bar-mu link-homotopy
  invariants of length n+1.  Second, we describe the structure of the
  Brunnian part of the degree 2n-graded quotient of the
  Goussarov--Vassiliev filtration for (n+1)-component links.
\end{asciiabstract}

\begin{htmlabstract}
The purpose of the paper is twofold.  First, we give a
short proof using the Kontsevich integral for the fact that the
restriction of an invariant of degree 2n to (n+1)&ndash;component
Brunnian links can be expressed as a quadratic form on the Milnor
&micro;-bar link-homotopy invariants of length n+1.  Second, we
describe the structure of the Brunnian part of the degree&ndash;2n
graded quotient of the Goussarov&ndash;Vassiliev filtration for
(n+1)&ndash;component links.
\end{htmlabstract}

\begin{webabstract}
The purpose of the paper is twofold.  First, we give a
short proof using the Kontsevich integral for the fact that the
restriction of an invariant of degree $2n$ to $(n{+}1)$--component
Brunnian links can be expressed as a quadratic form on the Milnor
$\bar\mu$ link-homotopy invariants of length $n+1$.  Second, we
describe the structure of the Brunnian part of the degree--$2n$
graded quotient of the Goussarov--Vassiliev filtration for
$(n{+}1)$--component links.
\end{webabstract}

\begin{abstract}
  The purpose of the paper is twofold.  First, we give a
  short proof using the Kontsevich integral for the fact that the
  restriction of an invariant of degree $2n$ to $(n{+}1)$--component
  Brunnian links can be expressed as a quadratic form on the Milnor
  $\bar\mu$ link-homotopy invariants of length $n+1$.  Second, we
  describe the structure of the Brunnian part of the degree--$2n$
  graded quotient of the Goussarov--Vassiliev filtration for
  $(n{+}1)$--component links.
\end{abstract}

\maketitle

\section{Introduction}
\label{sec:introduction}

The notion of Goussarov--Vassiliev finite type link invariants provides a
unifying viewpoint on the various quantum link invariants; see
Bar-Natan \cite{BN}, Goussarov \cite{Gusarov:91,Gusarov:94},
Vassiliev \cite{Vassiliev}.
The definition mainly relies on a descending filtration
\begin{equation}
  \label{e7}
  \modZ \modL (m)=J_0(m)\supset J_1(m)\supset \cdots
\end{equation}
of the free abelian group $\modZ \modL (m)$ generated by the set $\modL(m)$ of isotopy classes of $m$--component,
oriented, ordered links in $S^3$.  Here each $J_n(m)$ is generated by alternating sums of
links over $n$ independent crossing changes.  
A homomorphism from $\modZ \modL (m)$ to an abelian group is said to
be an invariant of degree $\le n$ if it vanishes on $J_{n+1}(m)$.
The Kontsevich integral \cite{Kontsevich} is a universal rational-valued finite type link invariant.

A link in $S^3$ is said to be {\em Brunnian\/} if every proper sublink
of it is an unlink.  The first author \cite{Hb} proved that an
$(n{+}1)$--component Brunnian link cannot be distinguished from the
unlink by any invariant of degree $<2n$ with values in any abelian
group.  In \cite{HM1} the authors study the behavior of the
invariants of degree $2n$ on $(n{+}1)$--component Brunnian links.

The purpose of this paper is twofold.  First, we give an alternative, shorter
proof using the Kontsevich integral of \fullref{r6},
originally proved in \cite{HM1}, which relates the degree--$2n$
invariants for $(n{+}1)$--component links to Milnor $\bar\mu$ link-homotopy
invariants.  Second, we describe the structure of the {\em Brunnian 
part} of the degree--$2n$ graded quotient of the Goussarov--Vassiliev
filtration for $(n{+}1)$--component links.

For any $(n{+}1)$--component Brunnian link $L$, the Milnor $\bar\mu$
link-homotopy invariants $\bar\mu_{i_1,\ldots,i_{n+1}}(L)\in\modZ$ for
nonrepeating indices $i_1,\ldots,i_{n+1}\in\{1,2,\ldots,n+1\}$ are 
well-defined and give a complete set of link-homotopy invariants
\cite{M}.  It is also known that there are exactly $(n-1)!$
independent invariants
\begin{equation*}
  \bar\mu _\sigma (L):= \bar\mu _{\sigma (1),\sigma (2),\ldots ,\sigma
  (n-1),n,n+1}(L)\in \modZ\quad \text{for $\sigma\in S_{n-1}$},
\end{equation*}
where $S_{n-1}$ denotes the symmetric group on
$\{1,2,\ldots,n-1\}$.

\begin{theorem}[{\cite[Theorem 1.2]{HM1}}, stated slightly differently]
  \label{r6}
  Let $n\ge 2$ and 
  let $f\colon\modZ\modL(n+1)\to\modZ$ be any link invariant of degree
  $2n$.  Then there exist integers $f_{\sigma ,\sigma '}$ for $\sigma
  ,\sigma '\in S_{n-1}$ such that for any $(n{+}1)$--component Brunnian
  link $L$, we have
  \begin{equation}
    \label{e8}
    f(L)-f(U)= \half\sumss f_{\sigma ,\sigma '}\bar\mu _\sigma (L)\bar\mu _{\sigma '}(L),
  \end{equation}
  where $U$ is the $(n{+}1)$--component unlink. 
  Moreover, we have $f_{\sigma ,\sigma '}=f_{\sigma ',\sigma}$ and $f_{\sigma ,\sigma}\in 2\mathbb{Z}$ for all 
  $\sigma,\sigma '\in S_{n-1}$.  
\end{theorem}

\fullref{r6} is proved in \fullref{sec:an-alternative-proof}.
The coefficients $f_{\sigma ,\sigma '}$ are defined there using the
weight system of $f$.  
It follows from the arguments of \fullref{sec:brunnian-part-brj} 
that we have $f_{\sigma ,\sigma '}=f_{\eta\sigma ,\eta\sigma '}$ for every $\eta\in S_{n-1}$.  
Thus there are integers $g_{\tau}$ for $\tau \in S_{n-1}$ 
satisfying $g_{\tau}=g_{\tau^{-1}}$ and $g_1\in 2\mathbb{Z}$, such that   
$f_{\sigma ,\sigma '}=g_{\sigma^{-1}\sigma '}$.  
Further, $f$ does not depend on the order of components for Brunnian links.  

In \fullref{sec:brunnian-part-brj} we give, in theory, a
necessary and sufficient condition to extend a quadratic form of Milnor
invariants to a link invariant of degree $\le 2n$.
More precisely, we study the {\em Brunnian
part} $\BrJ$ of the graded quotient $\bJ_{2n}(n+1)=J_{2n}(n+1)/J_{2n+1}(n+1)$, which is defined as the subgroup
generated by the elements $[L-U]_{J_{2n+1}}$, where $L$ is Brunnian.
We construct a homomorphism
\begin{equation*}
  h_n\zzzcolon \modA ^c_{n-1}(\emptyset)\rightarrow \BrJ,
\end{equation*}
where $\modA ^c_{n-1}(\emptyset)$ is a $\modZ $--module of connected trivalent
diagrams with $2n-2$ vertices.  We show that $h_n$ is surjective for
$n\ge 3$ and is an isomorphism over $\modQ $ for $n\ge 2$.

A few remarks are in order.

\begin{remark}
  \label{r1}
  The original proof of \fullref{r6} involved calculus of
  claspers.  It provided, as a byproduct, a number of results on
  Brunnian (string) links and the link-homotopy relation.  The proof
  given here is shorter than the original one and may also look
  simpler depending on the reader's point of view.
\end{remark}

\begin{remark}
  \label{r3}
  Recall that Milnor invariants of length $n+1$ for string links are
  Goussa\-rov--Vassiliev invariants of degree $\le n$; see Bar-Natan \cite{BN2} and Lin \cite{Lin}.
  For links, Milnor invariants are not well-defined in general.
  However, \fullref{r6} shows that a quadratic expression in
  Milnor invariants, which is well-defined at least for
  $(n{+}1)$--component Brunnian links, may extend to a link invariant of
  degree $\le 2n$.

  For example, the case $n=2$ of \fullref{r6}
  states that any degree--$4$ invariant of $3$--component links
  restricts for Brunnian links to a scalar multiple of
  the square of the triple Milnor invariant $\mu_{123}(L)$.  Recall
  that $\mu_{123}(L)^2$ for $3$--component algebraically split links
  extends to the degree--$4$ coefficient of the Alexander--Conway
  polynomial, which is a degree--$4$ invariant; see Cochran \cite{Cochran}.
\end{remark}

\begin{remark}
  \label{r2}
  Recall that trivalent diagrams appear in the study of Ohtsuki finite
  type invariants of integral homology spheres; see Ohtsuki \cite{O}, Garoufalidis and Ohtsuki \cite{GO} and Le \cite{Le}.  The
  relation between \fullref{r10} and Ohtsuki finite type
  invariants is discussed in \cite{jb}.
\end{remark}

\begin{acknowledgments}
  The authors wish to thank Thang Le and Christine Lescop for helpful comments and conversations.
  The first author is partially supported by the Japan Society
  for the Promotion of Science, Grant-in-Aid for Young Scientists (B),
  16740033.  The second author is supported by a Postdoctoral
  Fellowship and a Grant-in-Aid for Scientific Research of the Japan
  Society for the Promotion of Science.
\end{acknowledgments}

\section[Proof of Theorem 1.1 using the Kontsevich integral]{Proof of \fullref{r6} using the Kontsevich integral}
\label{sec:an-alternative-proof}
In this section, we give a proof of \fullref{r6} using the Kontsevich integral.
For an introduction to the Kontsevich integral; see for example 
Bar-Natan \cite{BN}, Chmutov and Duzhin, \cite{CD}, Lescop \cite{lescop} and Ohtsuki \cite{O2}.

\subsection{Spaces of diagrams}
\label{sec:diagrams}

We recall the definition of the spaces of diagrams in which the Kontsevich integral
takes its values.

A {\em unitrivalent diagram\/} (also called {\em Feynman diagram\/} or
{\em Jacobi diagram\/}) is a finite graph $\Gamma $ with univalent and
trivalent vertices such that each trivalent vertex is equipped with a
cyclic order of the three incident edges.
Throughout the paper, we assume that every component of a unitrivalent
diagram has at least one univalent vertex.  The {\em degree\/} of
$\Gamma $ is half the number of vertices in $\Gamma $.

\begin{figure}[ht!]
\centering
\begin{picture}(0,0)%
\includegraphics{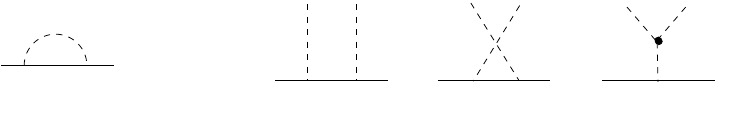}%
\end{picture}%
\setlength{\unitlength}{1973sp}%
\begingroup\makeatletter\ifx\SetFigFont\undefined%
\gdef\SetFigFont#1#2#3#4#5{%
  \reset@font\fontsize{#1}{#2pt}%
  \fontfamily{#3}\fontseries{#4}\fontshape{#5}%
  \selectfont}%
\fi\endgroup%
\begin{picture}(6986,1251)(1802,-4232)
\put(3064,-3397){\makebox(0,0)[lb]{\smash{\SetFigFont{6}{7.2}{\rmdefault}{\mddefault}{\updefault}{\color[rgb]{0,0,0}$=0$}%
}}}
\put(2202,-4184){\makebox(0,0)[lb]{\smash{\SetFigFont{6}{7.2}{\rmdefault}{\mddefault}{\updefault}{\color[rgb]{0,0,0}FI}%
}}}
\put(5590,-3491){\makebox(0,0)[lb]{\smash{\SetFigFont{6}{7.2}{\rmdefault}{\mddefault}{\updefault}{\color[rgb]{0,0,0}$-$}%
}}}
\put(7217,-3491){\makebox(0,0)[lb]{\smash{\SetFigFont{6}{7.2}{\rmdefault}{\mddefault}{\updefault}{\color[rgb]{0,0,0}$-$}%
}}}
\put(8788,-3491){\makebox(0,0)[lb]{\smash{\SetFigFont{6}{7.2}{\rmdefault}{\mddefault}{\updefault}{\color[rgb]{0,0,0}$=0$}%
}}}
\put(6255,-4232){\makebox(0,0)[lb]{\smash{\SetFigFont{6}{7.2}{\rmdefault}{\mddefault}{\updefault}{\color[rgb]{0,0,0}STU}%
}}}
\end{picture}\vspace{-2mm}
\caption{}
\label{F07graph}
\end{figure}

A unitrivalent diagram $\Gamma$ is {\em on a $1$--manifold $X$\/} if the univalent
vertices of $\Gamma$ are identified with distinct points on $X$ up to isotopy.
Let $\modA _k(X)$ denote the
$\modZ $--module generated by diagrams on $X$ of degree $k$, subject to the
{\em FI relations\/} and the {\em STU relations\/}; see \fullref{F07graph}.

\begin{figure}[ht!]
\centering
\begin{picture}(0,0)%
\includegraphics{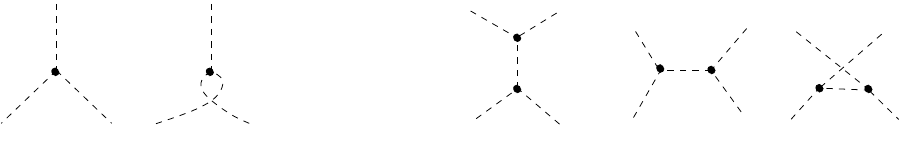}%
\end{picture}%
\setlength{\unitlength}{1973sp}%
\begingroup\makeatletter\ifx\SetFigFont\undefined%
\gdef\SetFigFont#1#2#3#4#5{%
  \reset@font\fontsize{#1}{#2pt}%
  \fontfamily{#3}\fontseries{#4}\fontshape{#5}%
  \selectfont}%
\fi\endgroup%
\begin{picture}(8734,1544)(2121,-6159)
\put(9353,-5238){\makebox(0,0)[lb]{\smash{\SetFigFont{6}{7.2}{\rmdefault}{\mddefault}{\updefault}{\color[rgb]{0,0,0}$+$}%
}}}
\put(10855,-5223){\makebox(0,0)[lb]{\smash{\SetFigFont{6}{7.2}{\rmdefault}{\mddefault}{\updefault}{\color[rgb]{0,0,0}$=0$}%
}}}
\put(7645,-5238){\makebox(0,0)[lb]{\smash{\SetFigFont{6}{7.2}{\rmdefault}{\mddefault}{\updefault}{\color[rgb]{0,0,0}$-$}%
}}}
\put(3212,-5238){\makebox(0,0)[lb]{\smash{\SetFigFont{6}{7.2}{\rmdefault}{\mddefault}{\updefault}{\color[rgb]{0,0,0}$+$}%
}}}
\put(4624,-5223){\makebox(0,0)[lb]{\smash{\SetFigFont{6}{7.2}{\rmdefault}{\mddefault}{\updefault}{\color[rgb]{0,0,0}$=0$}%
}}}
\put(8647,-6128){\makebox(0,0)[lb]{\smash{\SetFigFont{6}{7.2}{\rmdefault}{\mddefault}{\updefault}{\color[rgb]{0,0,0}IHX}%
}}}
\put(3242,-6159){\makebox(0,0)[lb]{\smash{\SetFigFont{6}{7.2}{\rmdefault}{\mddefault}{\updefault}{\color[rgb]{0,0,0}AS}%
}}}
\end{picture}\vspace{-2mm}
\caption{}
\label{F09}
\end{figure}

As is well-known, the {\em AS relations\/} and the {\em IHX relations\/}
depicted in \fullref{F09} are valid in $\modA _k(X)$.

If $X$ is the disjoint union of $m$ copies of $S^1$, then $\modA_k(X)$
is denoted also by $\modA _k(m)$.

An {\em open unitrivalent diagram\/} for $m$ colors $\{1,\ldots ,m\}$ is
a unitrivalent graph such that each univalent vertex is labeled by
an element of $\{1,\ldots ,m\}$.
Let $\modB _{k}(m)$ denote the $\modQ $--vector space generated by open unitrivalent diagrams of
degree $k$ for $m$ colors $\{1,\ldots ,m\}$, modulo the AS and the IHX relations.

Let $\sqcup^{m}I$ denote the disjoint union of $m$ copies of the unit
interval $I$.  There is a standard isomorphism ($m\ge 1$)
  \begin{equation*}
    \chi _{k}\zzzcolon \modB _{k}(m)\simeqto \modA _{k}(\sqcup^{m}I)\otimes \modQ ,
  \end{equation*}
called the Poincar\'e--Birkhoff--Witt isomorphism; see Bar-Natan \cite{BN}.

For $m=1$, $\modA _{k}(I)$ is isomorphic to $\modA _{k}(1)$.
For $m>1$, let $\modB ^l_{k}(m)$ denote the quotient of $\modB
_{k}(m)$ by all {\em link relations\/}
(see Section 5.2 of Bar-Natan, Garoufalidis, Rozansky and Thurston \cite{BNGRT}).  The isomorphism
$\chi _{k}$
induces a well-defined isomorphism
  \begin{equation*}
    \chi _{k}\zzzcolon \modB ^l_{k}(m)\simeqto \modA _{k}(m)\otimes \modQ.
  \end{equation*}

\subsection[Proof of Theorem 1.1]{Proof of \fullref{r6}}
\label{sec:newproof}

Suppose $L$ is an $(n{+}1)$--component Brunnian link in $S^3$.
By \cite[Proposition 12]{Hb}, there is an $(n{+}1)$--component Brunnian string link $T$
whose closure is $L$.  (A string link is said to be {\em Brunnian\/} if
every proper subtangle of it is a trivial string link.)
Consider the Kontsevich integral
$Z(T)$ of $T$, which lives in the completed $\modQ$--vector space
$\modA (\sqcup^{n+1}I)=\prod_{k\ge0} \modA _k(\sqcup^{n+1}I)\otimes\modQ$.
Recall that $\modA (\sqcup^{n+1}I)$ has an algebra
structure with multiplication given by the `stacking product' and with
unit $1$ given by the empty unitrivalent diagram.

Let $p\zzzcolon \modA (\sqcup^{n+1}I)\simeq\modB (n+1)$ denote the inverse of the
Poincar\'e--Birkhoff--Witt isomorphism, where $\modB (n+1)$ is the completed
$\modQ $--vector space $\prod_{k\ge0} \modB _k(n+1)$.
Since $p(Z(T))\in \modB (n+1)$ is grouplike (see Le, Murakami and Ohtsuki \cite{LMO}), we have
\begin{equation}
  \label{e4}
    p(Z(T)) =\exp_{\sqcup}(P),\quad P\in \modP (n+1),
\end{equation}
where $\modP (n+1)$ denotes the primitive part of $\modB (n+1)$, generated by
connected diagrams and where $\exp_{\sqcup}$ denotes the exponential with
respect to disjoint union of diagrams.

For each $i=1,\ldots ,n+1$, consider the operation $\epsilon _i$ of omitting the
$i$--th string.  At the level of the string link we have $\epsilon _i(T)=\modone _n$,
since $T$ is Brunnian.  Hence we have $\epsilon _i(P)=0$.  It follows that
$P$ can be expressed as an infinite $\modQ $--linear combination of
connected diagrams, each having at least one univalent vertex of color
$i$.  Since we have this property for $i=1,\ldots ,n+1$, we can deduce that
$P$ is an infinite $\modQ $--linear combination of connected diagrams, each
involving all the colors.

\begin{figure}[ht!]
\centering
\begin{picture}(0,0)%
\includegraphics{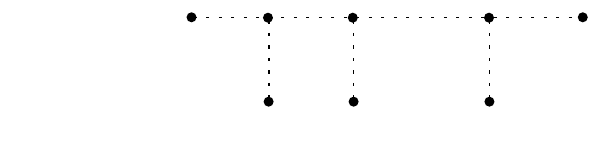}%
\end{picture}%
\setlength{\unitlength}{3947sp}%
\begingroup\makeatletter\ifx\SetFigFont\undefined%
\gdef\SetFigFont#1#2#3#4#5{%
  \reset@font\fontsize{#1}{#2pt}%
  \fontfamily{#3}\fontseries{#4}\fontshape{#5}%
  \selectfont}%
\fi\endgroup%
\begin{picture}(2890,695)(977,-386)
\put(977,-105){\makebox(0,0)[lb]{\smash{\SetFigFont{10}{12.0}{\familydefault}{\mddefault}{\updefault}{\color[rgb]{0,0,0}$t_\sigma=$}%
}}}
\put(2111,-341){\makebox(0,0)[lb]{\smash{\SetFigFont{10}{12.0}{\familydefault}{\mddefault}{\updefault}{\color[rgb]{0,0,0}$\sigma(1)$}%
}}}
\put(2563,-331){\makebox(0,0)[lb]{\smash{\SetFigFont{10}{12.0}{\familydefault}{\mddefault}{\updefault}{\color[rgb]{0,0,0}$\sigma(2)$}%
}}}
\put(3237,-338){\makebox(0,0)[lb]{\smash{\SetFigFont{10}{12.0}{\familydefault}{\mddefault}{\updefault}{\color[rgb]{0,0,0}$\sigma(n-1)$}%
}}}
\put(2833,-28){\makebox(0,0)[lb]{\smash{\SetFigFont{10}{12.0}{\familydefault}{\mddefault}{\updefault}{\color[rgb]{0,0,0}$\dots$}%
}}}
\put(1444,201){\makebox(0,0)[lb]{\smash{\SetFigFont{10}{12.0}{\familydefault}{\mddefault}{\updefault}{\color[rgb]{0,0,0}$n+1$}%
}}}
\put(3867,201){\makebox(0,0)[lb]{\smash{\SetFigFont{10}{12.0}{\familydefault}{\mddefault}{\updefault}{\color[rgb]{0,0,0}$n$}%
}}}
\end{picture}
\label{F32new}\vspace{-2mm}
\caption{}
\end{figure}

By an easy counting argument, we see that $P$ can be expressed as
\begin{equation*}
  P=P_n + P_{n+1}+ P_{n+2}+\cdots 
\end{equation*}
where $P_n$ is a linear combination of trees of degree $n$, and $P_k$
for $k>n$ is a linear combination of connected diagrams of degree $k$.
Moreover, we can use the AS and IHX relation to express $P_n$ as a
linear combination of trees $t_\sigma $ as depicted in \fullref{F32new},
$\sigma\in S_{n-1}$.
 By
Habegger and Masbaum's result \cite{HM}, $P_n$ is thus a $\modZ $--linear
combination of trees which corresponds to the Milnor invariants of
length $n+1$ of the form $\bar\mu _{\sigma}:=\bar\mu _{\sigma (1),\sigma (2),\ldots ,\sigma (n-1),n,n+1}$, for
$\sigma\in S_{n-1}$.
More precisely, from \cite[Theorem 6.1 and Proposition 7.1]{HM} and properties of
the Milnor $\bar\mu$ link-homotopy invariants \cite{M} it follows that
\begin{equation} \label{mu}
  P_n=\sum_{\sigma\in S_{n-1}} \bar\mu _{\sigma}(L)t_{\sigma}.
\end{equation}
By \eqref{e4}, we have
\begin{equation*}
  \begin{split}
    p(Z(T))
    &= \exp_{\sqcup}(P_n+P_{n+1}+\cdots +P_{2n}+\cdots )\\
    &= 1+P_n+P_{n+1}+\cdots +P_{2n-1}+(P_{2n}+\half P_n^2)+O_{>2n},
  \end{split}
\end{equation*}
where $O_{>2n}$ denotes a sum of terms of degree greater than $2n$.
Hence we have
\begin{equation}
  \label{e1}
  Z(T)
  = 1+P'_n+P'_{n+1}+\cdots +P'_{2n-1}+P'_{2n}+
  p^{-1}(\half P_n^2)+O_{>2n},
\end{equation}
where $P'_k=p^{-1}(P_k)$ for $k=n,\ldots ,2n$.

By \cite[Theorem 6.1]{LM}, the Kontsevich integral $Z(L)$ of $L$ can be obtained from $Z(T)$ by first multiplying
by the $(n{+}1)$--parallel $\Delta^n(Z(\circlearrowleft))$ of the Kontsevich integral of the unknot
$Z(\circlearrowleft) \in \modA (I)$ and then mapping into the space $\modA (n+1)$ by closing the
strings.  Thus we have
\begin{equation}
  \label{e5}
  Z(L)-Z(U)=\pi \left( \Delta^n(Z(\circlearrowleft))(Z(T)-1) \right),
\end{equation}
where $U$ denotes the $(n{+}1)$--component unlink in the $3$--sphere $S^3$ and where
the projection $\pi \zzzcolon \modA (\sqcup^{n+1}I)\rightarrow \modA (n+1)$ is 
induced by the closure operation.

As is well-known \cite{BNGRT2,thurston,BLT},
\begin{equation}
  \label{e2}
  Z(\circlearrowleft)= 1+ W_2+ W_4+W_6+ \cdots
\end{equation}
where $W_{2j}$ is a linear combination of unitrivalent diagram of
degree $2j$ consisting of finitely many {\em wheels\/}, ie,
connected unitrivalent diagrams of Euler number $1$.

By \eqref{e1}, \eqref{e5} and \eqref{e2}, we have
\begin{equation}
  \label{e6}
  Z(L)-Z(U)
  =
  \half\pi p^{-1}(P_n^2)
  +\sum_{i=n}^{2n}\pi(P'_i)
  +\sum_{i,j}\pi(P'_i\Delta^n(W_{2j})) + O_{>2n},
\end{equation}
where the second sum is over $i,j$ with $n\le i\le 2n$, $2\le 2j\le n$,
$i+2j\le 2n$.

One can verify that the two sums in the right-hand side of \eqref{e6}
vanish, using the fact that  $\pi$ kills each diagram with at least one string
having exactly one univalent vertex.  Thus we have
\begin{equation} \label{Z}
    Z(L)-Z(U) = \half\pi p^{-1}(P_n^2)+O_{>2n}.
\end{equation}
It follows that the first nontrivial, $\modZ $--valued, finite type invariants of $(n{+}1)$--compo\-nent
Brunnian links are of degree $2n$.
Moreover, \eqref{mu} and \eqref{Z} imply that
\begin{equation*}
  Z(L)-Z(U)=   \half\sumss \bar\mu _{\sigma}(L)\bar\mu _{\sigma'}(L) \pi p^{-1}(t_{\sigma}t_{\sigma'}) +
      O_{>2n}.
\end{equation*}
Let $f$ be a $\modZ $--valued link invariant of degree $2n$, and
denote by $W_f$ the weight system of $f$.
To complete the proof of \fullref{r6}, we obtain \eqref{e8} by setting
\begin{equation*}
  f_{\sigma ,\sigma '} =
  W_f\left(\pi p^{-1}(t_{\sigma}t_{\sigma'})\right).
\end{equation*}

\begin{remark}
\label{rem:coeff}
Obviously, we have $f_{\sigma,\sigma'}=f_{\sigma',\sigma}$ for all
$\sigma, \sigma'\in S_{n-1}$.  Observe that the identity depicted in \fullref{F28} can be derived from the STU relation.  
This identity implies that, for any $\sigma\in S_{n-1}$, we have $\pi
p^{-1}(t_{\sigma}^2)=2d_{\sigma}$ for some $d_{\sigma}\in\modA
_{2n}(n+1)$.  
\begin{figure}[ht!]
\centering
\begin{picture}(0,0)%
\includegraphics{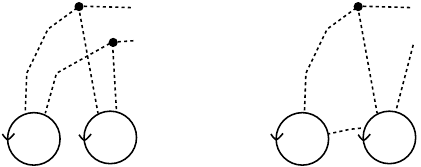}%
\end{picture}%
\setlength{\unitlength}{3947sp}%
\begingroup\makeatletter\ifx\SetFigFont\undefined%
\gdef\SetFigFont#1#2#3#4#5{%
  \reset@font\fontsize{#1}{#2pt}%
  \fontfamily{#3}\fontseries{#4}\fontshape{#5}%
  \selectfont}%
\fi\endgroup%
\begin{picture}(2004,795)(762,-1399)
\put(1504,-1061){\makebox(0,0)[lb]{\smash{\SetFigFont{12}{14.4}{\familydefault}{\mddefault}{\updefault}{\color[rgb]{0,0,0}$=\;-2$}%
}}}
\end{picture}
\caption{Hint: use the fact that, thanks to the STU relation, a circle component inserted in an edge 
attached to a trivalent vertex $v$ can be moved to another edge attached to $v$ (see \fullref{F23}). } 
\label{F28}
\end{figure}
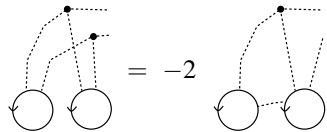
Hence $f_{\sigma ,\sigma}=W_f(2d_{\sigma})$ is even.
\end{remark}

\begin{remark}
As recalled in the introduction, it is known that no
Goussarov--Vassiliev invariant of degree $<2n$ with values in any abelian group can distinguish
$L$ from $U$ \cite{Hb}.
The above proof of \fullref{r6} recovers, using the Kontsevich integral, a weaker version (namely,
a rational version) of this result.
\end{remark}

\section[The Brunnian part of the 2n-th graded part of the Goussarov--Vassiliev filtration]{The Brunnian part of $\bJ_{2n}(n+1)$}
\label{sec:brunnian-part-brj}

The purpose of this section is to give an almost complete description
of the structure of the Brunnian part $\BrJ$ of the $2n$--th graded
quotient $\bJ_{2n}(n+1)$ of the Goussarov--Vassiliev filtration for 
$(n{+}1)$--component links, using the $\modZ $--module $\modA ^c_{n-1}(\emptyset)$
of connected trivalent diagrams.

\subsection{Claspers and Goussarov--Vassiliev filtration}
\label{sec:graph-schem-gouss}

The Goussarov--Vassiliev filtration \eqref{e7} is usually defined in
terms of singular links, where $J_k(m)\subset\modZ\modL(m)$ is generated by
the $k$--fold alternating sums of links determined by singular links with $k$
double points; see eg Bar-Natan \cite{BN}.
One can also redefine the filtration using claspers; see Habiro \cite{H} and Goussarov \cite{Goussarov01}.
Here we will freely use the definitions and conventions of \cite{H}.
Recall that a {\em graph scheme of degree $k$\/} for a link $L$ in $S^3$ is a set
$S=\{G_1,\ldots ,G_l\}$ of  $l$ disjoint (strict) graph claspers $G_1,\ldots,G_l$ for $L$ such that $\sum_{i=1}^l \deg G_i=k$.
We set
\begin{equation}
 \label{eq:alt}
  [L,S]=[L;G_1,\ldots,G_l] := \sum_{S'\subset S} (-1)^{|S|-|S'|} L_{\bigcup S'} \in \modZ \modL (m),
\end{equation}
where the sum runs over all subsets $S'$ of $S$, $|S'|$ denote the
number of elements of $S'$, and $L_{\bigcup S'}$ denotes the result
from $L$ of surgery along the union of elements in $S'$.
It is known that $J_k(m)$ is generated by the elements of the form
$[L,S]$, where $L\in \modL (m)$ and $S$ is a graph scheme for $L$ of degree $k$ \cite{H}.

\subsection{The Brunnian part of the Goussarov--Vassiliev filtration}
\label{brunnianpart}

If $L$ is an $(n{+}1)$--component Brunnian link in $S^3$, then, by
\cite[Theorem 3]{Hb}, we have $L-U\in J_{2n}(n+1)$.  Define $\BrJ$ to
be the $\modZ $--submodule of $\bJ_{2n}(n+1)$ generated by the elements
$[L-U]_{J_{2n+1}}$, where $L$ is an $(n{+}1)$--component Brunnian link.
We call $\BrJ$ the {\em Brunnian part\/} of $\bJ_{2n}(n+1)$.

\subsection{Trivalent diagrams}
\label{sec:trivalent-diagrams}

A \emph{trivalent diagram\/\/} is a finite trivalent graph $\Gamma$ such
that each vertex is equipped with a cyclic order on the three incident
edges.  The \emph{degree\/\/} of $\Gamma$ is half the number of its vertices.

For $k\ge 0$, let $\mathcal{A}_k(\emptyset)$ denote the $\mathbb{Z\/}$--module generated by trivalent diagrams of degree
$k$, subject to the AS and IHX relations (see \fullref{F09}).

Let $\modA ^c_n(\emptyset)$ denote the $\modZ $--submodule of $\modA _n(\emptyset)$
generated by connected trivalent diagrams.

\subsection[The circle-insertion map]{The circle-insertion map $g_n$}
\label{sec:circle-insertion-map}

\begin{figure}[ht!]
\centering
\begin{picture}(0,0)%
\includegraphics{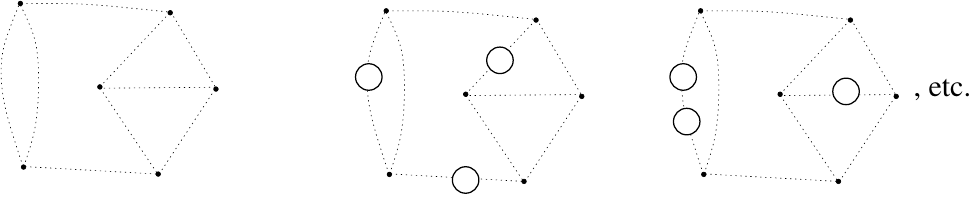}%
\end{picture}%
\setlength{\unitlength}{1579sp}%
\begingroup\makeatletter\ifx\SetFigFont\undefined%
\gdef\SetFigFont#1#2#3#4#5{%
  \reset@font\fontsize{#1}{#2pt}%
  \fontfamily{#3}\fontseries{#4}\fontshape{#5}%
  \selectfont}%
\fi\endgroup%
\begin{picture}(11699,2658)(1790,-4755)
\put(9606,-2850){\makebox(0,0)[lb]{\smash{\SetFigFont{7}{8.4}{\rmdefault}{\mddefault}{\updefault}{\color[rgb]{0,0,0}\adjustlabel<-1pt,-2pt>{$2$}}%
}}}
\put(11639,-3517){\makebox(0,0)[lb]{\smash{\SetFigFont{7}{8.4}{\rmdefault}{\mddefault}{\updefault}{\color[rgb]{0,0,0}\adjustlabel<-2pt,-2pt>{$3$}}%
}}}
\put(9590,-3683){\makebox(0,0)[lb]{\smash{\SetFigFont{7}{8.4}{\rmdefault}{\mddefault}{\updefault}{\color[rgb]{0,0,0}$1$}%
}}}
\put(4670,-3182){\makebox(0,0)[lb]{\smash{\SetFigFont{10}{12.0}{\familydefault}{\mddefault}{\updefault}{\color[rgb]{0,0,0}$\overset{g}{\mapsto}$}%
}}}
\put(5770,-2966){\makebox(0,0)[lb]{\smash{\SetFigFont{7}{8.4}{\rmdefault}{\mddefault}{\updefault}{\color[rgb]{0,0,0}$1$}%
}}}
\put(7253,-4683){\makebox(0,0)[lb]{\smash{\SetFigFont{7}{8.4}{\rmdefault}{\mddefault}{\updefault}{\color[rgb]{0,0,0}\adjustlabel<0pt,-2pt>{$2$}}%
}}}
\put(8003,-2900){\makebox(0,0)[lb]{\smash{\SetFigFont{7}{8.4}{\rmdefault}{\mddefault}{\updefault}{\color[rgb]{0,0,0}$3$}%
}}}
\put(9076,-3286){\makebox(0,0)[lb]{\smash{\SetFigFont{10}{12.0}{\familydefault}{\mddefault}{\updefault}{\color[rgb]{0,0,0}$=$}%
}}}
\end{picture}
\label{F23}
\caption{}
\end{figure}

For $n\ge 2$, denote the homomorphism which maps each trivalent diagram $\Gamma $ to the
result of inserting $n+1$ ordered copies of $S^1$ in its edges by
\begin{equation*}
  g_n\zzzcolon \modA ^c_{n-1}(\emptyset) \rightarrow  \modA _{2n\/}(n+1).
\end{equation*}
 See
\fullref{F23}. 
This map is well-defined thanks to the
STU relation, since $\Gamma $ is connected.  

We need the following result.

\begin{proposition}
  \label{r13}
  For $n\ge 2$, the map
  \begin{equation*}
    g_n\otimes \modQ \zzzcolon \modA ^c_{n-1}(\emptyset)\otimes \modQ \rightarrow \modA _{2n\/}(n+1)\otimes \modQ
  \end{equation*}
  is injective.
\end{proposition}

\begin{proof}
  Define a homomorphism
  \begin{equation*}
    P\zzzcolon \modB ^l_{2n}(n+1) \longrightarrow\modA _{n-1}(\emptyset)\otimes \modQ
  \end{equation*}
  as follows.  Let $\Gamma$ be a diagram in $\modB ^l_{2n}(n+1)$.  If there are
  exactly two univalent vertices colored by $i$ for each
  $i=1,\ldots ,n+1$, then we set $P(\Gamma )$ to be the trivalent diagram
  obtained from $\Gamma $ by joining each pair of univalent vertices of the
  same color.  Otherwise, set $P(\Gamma )=0$.  Let
  \begin{equation*}
    \pi \zzzcolon \modA _{n-1}(\emptyset)\otimes \modQ \longrightarrow \modA ^c_{n-1\/}(\emptyset)\otimes \modQ
  \end{equation*}
  denote the projection which maps each connected diagram into
  itself and maps nonconnected diagrams to $0$.  One can easily check
  that the composition
  \begin{equation*}
    \modA ^c_{n-1}(\emptyset)\otimes \modQ
    \negthinspace \xto{g_n\otimes \modQ }\negthinspace  \modA _{2n}(n+1)\otimes \modQ
    \negthinspace \xto{\chi _{2n}^{-1}}\negthinspace \modB ^l_{2n}(n+1)
    \negthinspace \xto{P}\negthinspace  \modA _{n-1}(\emptyset)\otimes \modQ
    \negthinspace \xto{\pi }\negthinspace  \modA ^c_{n-1}(\emptyset)\otimes \modQ
  \end{equation*}
  is the identity.  Hence $g_n\otimes \modQ $ is injective.
\end{proof}

\subsection[Structure of the Brunnian part]{Structure of $\BrJ$}
\label{sec:statement-result}

\begin{figure}[ht!]
\centering
\begin{picture}(0,0)%
\includegraphics{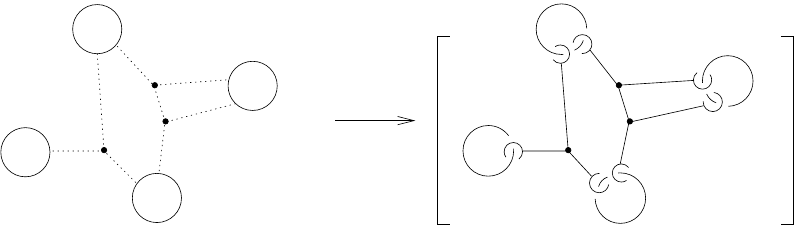}%
\end{picture}%
\setlength{\unitlength}{1579sp}%
\begingroup\makeatletter\ifx\SetFigFont\undefined%
\gdef\SetFigFont#1#2#3#4#5{%
  \reset@font\fontsize{#1}{#2pt}%
  \fontfamily{#3}\fontseries{#4}\fontshape{#5}%
  \selectfont}%
\fi
\endgroup%
\begin{picture}(9594,2797)(232,-2753)
\put(288,-2239){\makebox(0,0)[lb]{\smash{\SetFigFont{6}{7.2}{\rmdefault}{\mddefault}{\updefault}{\color[rgb]{0,0,0}\adjustlabel<0pt,-3pt>{$1$}}%
}}}
\put(1734,-350){\makebox(0,0)[lb]{\smash{\SetFigFont{6}{7.2}{\rmdefault}{\mddefault}{\updefault}{\color[rgb]{0,0,0}$2$}%
}}}
\put(2431,-2567){\makebox(0,0)[lb]{\smash{\SetFigFont{6}{7.2}{\rmdefault}{\mddefault}{\updefault}{\color[rgb]{0,0,0}$3$}%
}}}
\put(3188,-1458){\makebox(0,0)[lb]{\smash{\SetFigFont{6}{7.2}{\rmdefault}{\mddefault}{\updefault}{\color[rgb]{0,0,0}\adjustlabel<0pt,-3pt>{$4$}}%
}}}
\put(1167,-2050){\makebox(0,0)[lb]{\smash{\SetFigFont{6}{7.2}{\rmdefault}{\mddefault}{\updefault}{\color[rgb]{0,0,0}\adjustlabel<0pt,-2pt>{$\Gamma_1$}}%
}}}
\put(2119,-818){\makebox(0,0)[lb]{\smash{\SetFigFont{6}{7.2}{\rmdefault}{\mddefault}{\updefault}{\color[rgb]{0,0,0}$\Gamma_2$}%
}}}
\put(5867,-2239){\makebox(0,0)[lb]{\smash{\SetFigFont{6}{7.2}{\rmdefault}{\mddefault}{\updefault}{\color[rgb]{0,0,0}\adjustlabel<0pt,-2pt>{$1$}}%
}}}
\put(7328,-185){\makebox(0,0)[lb]{\smash{\SetFigFont{6}{7.2}{\rmdefault}{\mddefault}{\updefault}{\color[rgb]{0,0,0}$2$}%
}}}
\put(8011,-2559){\makebox(0,0)[lb]{\smash{\SetFigFont{6}{7.2}{\rmdefault}{\mddefault}{\updefault}{\color[rgb]{0,0,0}$3$}%
}}}
\put(6778,-2000){\makebox(0,0)[lb]{\smash{\SetFigFont{6}{7.2}{\rmdefault}{\mddefault}{\updefault}{\color[rgb]{0,0,0}\adjustlabel<-1pt,-4pt>{$G_1$}}%
}}}
\put(7739,-827){\makebox(0,0)[lb]{\smash{\SetFigFont{6}{7.2}{\rmdefault}{\mddefault}{\updefault}{\color[rgb]{0,0,0}\adjustlabel<-2pt,0pt>{$G_2$}}%
}}}
\put(8955,-1466){\makebox(0,0)[lb]{\smash{\SetFigFont{6}{7.2}{\rmdefault}{\mddefault}{\updefault}{\color[rgb]{0,0,0}\adjustlabel<0pt,-2pt>{$4$}}%
}}}
\put(4351,-1111){\makebox(0,0)[lb]{\smash{\SetFigFont{6}{7.2}{\rmdefault}{\mddefault}{\updefault}{\color[rgb]{0,0,0}\adjustlabel<3pt,-1pt>{$\xi_k$}}%
}}}
\put(9826,-2686){\makebox(0,0)[lb]{\smash{\SetFigFont{6}{7.2}{\rmdefault}{\mddefault}{\updefault}{\color[rgb]{0,0,0}alt,$J_{k+1}$}%
}}}
\end{picture}
\caption{Realization of unitrivalent diagrams and the map $\smash{\xi_k}$.
Here $[-]_{\text{alt},J_{k+1}}$ denotes the $J_{k+1}$--coset of the
alternating sum of the form \eqref{eq:alt} defined by the link and graph claspers.}  
\label{F00new}
\end{figure}
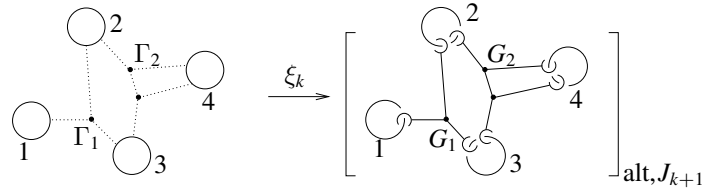

We need the well-known surjective surgery map
\begin{equation*}
  \xi _{k}\zzzcolon \modA _{k}(m) \rightarrow \bJ_{k}(m)
\end{equation*}
which maps each unitrivalent diagram
$\Gamma=\Gamma_1\cup\cdots\cup\Gamma_l$ to the coset of alternating
sum $[U;G_1,\ldots,G_l]$, where $G_1,\ldots,G_l$ are ``realizations''
of $\Gamma_1,\ldots,\Gamma_l$; see \fullref{F00new}.
This map coincides with the standard map which ``replaces chords with double points'' \cite[Section 8.2]{H}.  
$\xi _{k}\otimes\modQ$ is an isomorphism, with inverse given by the Kontsevich integral.  

For $n\ge 2$, let $h_n$ denote the composition
\begin{equation*}
  h_n\zzzcolon \modA ^c_{n-1}(\emptyset)\xto{g_n\/}\modA _{2n}(n+1)\xto{\xi _{2n}}\bJ_{2n}(n+1).
\end{equation*}

\begin{theorem}
  \label{r10}
  {\rm (1)}\qua For $n\ge 3$ we have $h_n(\modA ^c_{n-1}(\emptyset))=\BrJ$.  For
  $n=2$, $h_2(\modA ^c_1(\emptyset))$ is an index $2$ subgroup of
  $\Br(\bJ_4(3))$.

  {\rm (2)}\qua For $n\ge 2$, the $\modQ $--linear map
  \begin{equation}
    \label{e3}
    h_n\otimes \modQ \zzzcolon \modA ^c_{n-1}(\emptyset)\otimes \modQ \rightarrow \Br(\bJ_{2n\/}(n+1))\otimes \modQ
  \end{equation}
  is an isomorphism.
\end{theorem}

\begin{proof}
  (1)\qua Suppose $n\ge 3$.  First we show that
  \begin{equation}
    \label{e10}
    \BrJ\subset h_n(\modA ^c_{n-1}(\emptyset)).
  \end{equation}
  As outlined in \cite[Section 7]{HM1}, $\Br(\bar{J}_{2n}(n+1))$ is
  $\modZ $--spanned by
  \begin{align*}
    \half[U;T_\sigma ,\tilde T_\sigma ]\quad &\text{for $\sigma \in S_{n-1}$},\\
         [U;T_\sigma ,\tilde T_{\sigma '}]\quad &\text{for $\sigma ,\sigma '\in S_{n-1}$},
  \end{align*}
  where for all $\sigma,\sigma'$ in the symmetric group $S_{n-1}$, $T_{\sigma}$ is the simple degree--$n$ tree clasper for the
  $(n{+}1)$--component unlink $U$ depicted in \fullref{Tsigma}, and $\tilde{T}_{\sigma'}$ is obtained from $T_{\sigma'}$ by a small
  isotopy so that it is disjoint from $T_{\sigma}$.  
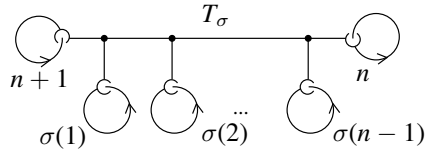
\begin{figure}[ht!]
\centering
  \begin{picture}(0,0)%
\includegraphics{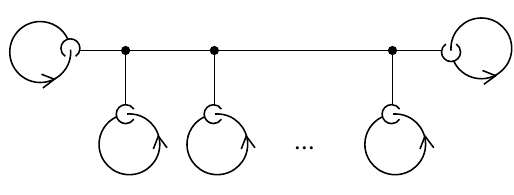}%
\end{picture}%
\setlength{\unitlength}{3947sp}%
\begingroup\makeatletter\ifx\SetFigFont\undefined%
\gdef\SetFigFont#1#2#3#4#5{%
  \reset@font\fontsize{#1}{#2pt}%
  \fontfamily{#3}\fontseries{#4}\fontshape{#5}%
  \selectfont}%
\fi\endgroup%
\begin{picture}(2465,927)(1,-406)
\put(1201,389){\makebox(0,0)[lb]{\smash{\SetFigFont{10}{12.0}{\rmdefault}{\mddefault}{\updefault}{\color[rgb]{0,0,0}$T_{\sigma}$}%
}}}
\put(2026,-361){\makebox(0,0)[lb]{\smash{\SetFigFont{10}{12.0}{\rmdefault}{\mddefault}{\updefault}{\color[rgb]{0,0,0}$\sigma(n-1)$}%
}}}
\put(2176, 14){\makebox(0,0)[lb]{\smash{\SetFigFont{10}{12.0}{\rmdefault}{\mddefault}{\updefault}{\color[rgb]{0,0,0}$n$}%
}}}
\put(  1, 14){\makebox(0,0)[lb]{\smash{\SetFigFont{10}{12.0}{\rmdefault}{\mddefault}{\updefault}{\color[rgb]{0,0,0}\adjustlabel<0pt,-2pt>{$n+1$}}%
}}}
\put(1201,-361){\makebox(0,0)[lb]{\smash{\SetFigFont{10}{12.0}{\rmdefault}{\mddefault}{\updefault}{\color[rgb]{0,0,0}$\sigma(2)$}%
}}}
\put(201,-361){\makebox(0,0)[lb]{\smash{\SetFigFont{10}{12.0}{\rmdefault}{\mddefault}{\updefault}{\color[rgb]{0,0,0}\adjustlabel<-1pt,-1pt>{$\sigma(1)$}}%
}}}
\end{picture}
\caption{The tree clasper $T_{\sigma}$}
\label{Tsigma}
\end{figure}
  ($\half[U;T_\sigma ,\tilde T_\sigma ]$ means an element $x\in\BrJ$
  such that $2x=[U;T_\sigma ,\tilde T_\sigma ]$.  Existence of such an
  element is shown in \cite{HM1}.)

  For $n\ge 3$, one can check using the AS and IHX relations that $\half[U;T_\sigma ,\tilde T_\sigma ]_{J_{2n+1}}$ is equal to:
  \begin{equation*} \def\xxab{$\half[U;T_\sigma ,\tilde T_\sigma
  ]_{J_{2n+1}}=\pm \Biggl[\hspace{53mm}\Biggr]_{J_{2n+1}}$}
  \def\xxb{$=\pm h_n\Biggl(\hspace{37mm}\Biggr)$}
\begin{picture}(0,0)%
\includegraphics{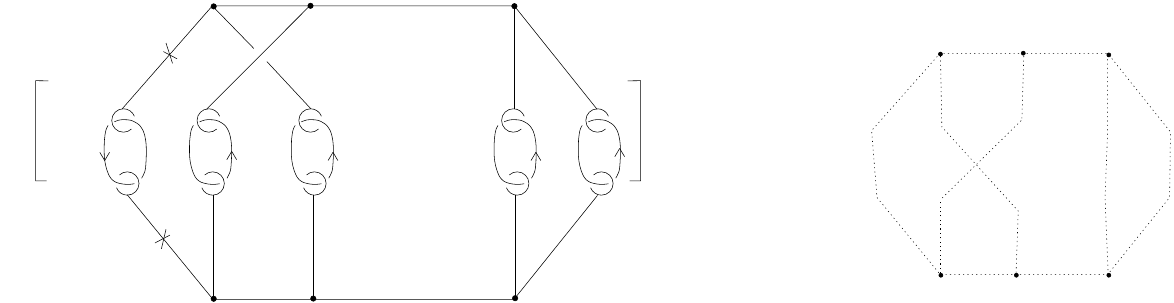}%
\end{picture}%
\setlength{\unitlength}{1381sp}%
\begingroup\makeatletter
\gdef\SetFigFont#1#2#3#4#5{%
  \reset@font\fontsize{#1}{#2pt}%
  \fontfamily{#3}\fontseries{#4}\fontshape{#5}%
  \selectfont}%
\endgroup%
\begin{picture}(16066,4108)(1932,-11975)
\put(16369,-9355){\makebox(0,0)[lb]{\smash{\SetFigFont{5}{6.0}{\familydefault}{\mddefault}{\updefault}{\color[rgb]{0,0,0}$\cdots$}%
}}}
\put(16243,-11088){\makebox(0,0)[lb]{\smash{\SetFigFont{5}{6.0}{\familydefault}{\mddefault}{\updefault}{\color[rgb]{0,0,0}$\cdots$}%
}}}
\put(11802,-9985){\makebox(0,0)[lb]{\smash{\SetFigFont{6}{7.2}{\familydefault}{\mddefault}{\updefault}{\color[rgb]{0,0,0}\adjustlabel<8pt,0pt>{\xxb}}%
}}}
\put(7382,-9775){\makebox(0,0)[lb]{\smash{\SetFigFont{6}{7.2}{\familydefault}{\mddefault}{\updefault}{\color[rgb]{0,0,0}$\sigma(n-1)$}%
}}}
\put(7335,-11260){\makebox(0,0)[lb]{\smash{\SetFigFont{6}{7.2}{\familydefault}{\mddefault}{\updefault}{\color[rgb]{0,0,0}$\cdots$}%
}}}
\put(10755,-10510){\makebox(0,0)[lb]{\smash{\SetFigFont{6}{7.2}{\familydefault}{\mddefault}{\updefault}{\color[rgb]{0,0,0}$J_{2n+1}$}%
}}}
\put(7293,-8961){\makebox(0,0)[lb]{\smash{\SetFigFont{6}{7.2}{\familydefault}{\mddefault}{\updefault}{\color[rgb]{0,0,0}$\cdots$}%
}}}
\put(1932,-9925){\makebox(0,0)[lb]{\smash{\SetFigFont{6}{7.2}{\familydefault}{\mddefault}{\updefault}{\color[rgb]{0,0,0}$\pm$}%
}}}
\put(6526,-10486){\makebox(0,0)[lb]{\smash{\SetFigFont{6}{7.2}{\familydefault}{\mddefault}{\updefault}{\color[rgb]{0,0,0}$\sigma(2)$}%
}}}
\put(9602,-9719){\makebox(0,0)[lb]{\smash{\SetFigFont{6}{7.2}{\familydefault}{\mddefault}{\updefault}{\color[rgb]{0,0,0}$n$}%
}}}
\put(5116,-10479){\makebox(0,0)[lb]{\smash{\SetFigFont{6}{7.2}{\familydefault}{\mddefault}{\updefault}{\color[rgb]{0,0,0}$\sigma(1)$}%
}}}
\put(2534,-9625){\makebox(0,0)[lb]{\smash{\SetFigFont{6}{7.2}{\familydefault}{\mddefault}{\updefault}{\color[rgb]{0,0,0}$n+1$}%
}}}
\end{picture}
   \end{equation*} For $n\ge 2$, $[U;T_\sigma
  ,\tilde T_{\sigma '}]_{J_{2n+1}}$ is equal to: \begin{equation}
  \label{e12} \def\xxa{$[U;T_\sigma ,\tilde T_{\sigma
  '}]_{J_{2n+1}}=\pm $} \def\xxb{$=\pm h_n\Biggl(\hspace{37mm}\Biggr)$}
\begin{picture}(0,0)%
\includegraphics{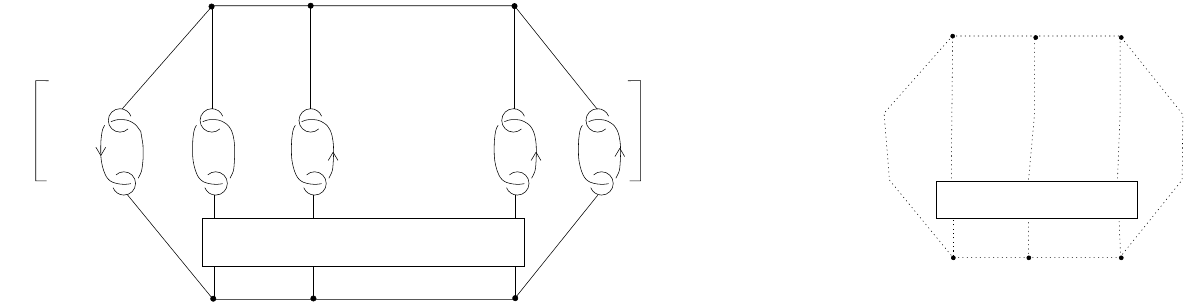}%
\end{picture}%
\setlength{\unitlength}{1381sp}%
\begingroup\makeatletter
\gdef\SetFigFont#1#2#3#4#5{%
  \reset@font\fontsize{#1}{#2pt}%
  \fontfamily{#3}\fontseries{#4}\fontshape{#5}%
  \selectfont}%
\endgroup%
\begin{picture}(16235,4108)(629,-16474)
\put(10586,-14246){\makebox(0,0)[lb]{\smash{\SetFigFont{6}{7.2}{\familydefault}{\mddefault}{\updefault}{\color[rgb]{0,0,0}\adjustlabel<8pt,0pt>{\xxb}}%
}}}
\put(15235,-13616){\makebox(0,0)[lb]{\smash{\SetFigFont{5}{6.0}{\familydefault}{\mddefault}{\updefault}{\color[rgb]{0,0,0}$\cdots$}%
}}}
\put(15088,-15714){\makebox(0,0)[lb]{\smash{\SetFigFont{5}{6.0}{\familydefault}{\mddefault}{\updefault}{\color[rgb]{0,0,0}$\cdots$}%
}}}
\put(6079,-14274){\makebox(0,0)[lb]{\smash{\SetFigFont{6}{7.2}{\familydefault}{\mddefault}{\updefault}{\color[rgb]{0,0,0}$\sigma(n-1)$}%
}}}
\put(5990,-13460){\makebox(0,0)[lb]{\smash{\SetFigFont{6}{7.2}{\familydefault}{\mddefault}{\updefault}{\color[rgb]{0,0,0}$\cdots$}%
}}}
\put(8299,-14218){\makebox(0,0)[lb]{\smash{\SetFigFont{6}{7.2}{\familydefault}{\mddefault}{\updefault}{\color[rgb]{0,0,0}$n$}%
}}}
\put(629,-14424){\makebox(0,0)[lb]{\smash{\SetFigFont{6}{7.2}{\familydefault}{\mddefault}{\updefault}{\color[rgb]{0,0,0}$\pm$}%
}}}
\put(5176,-14986){\makebox(0,0)[lb]{\smash{\SetFigFont{6}{7.2}{\familydefault}{\mddefault}{\updefault}{\color[rgb]{0,0,0}$\sigma(2)$}%
}}}
\put(3826,-14986){\makebox(0,0)[lb]{\smash{\SetFigFont{6}{7.2}{\familydefault}{\mddefault}{\updefault}{\color[rgb]{0,0,0}$\sigma(1)$}%
}}}
\put(1231,-14124){\makebox(0,0)[lb]{\smash{\SetFigFont{6}{7.2}{\familydefault}{\mddefault}{\updefault}{\color[rgb]{0,0,0}$n+1$}%
}}}
\put(9555,-14885){\makebox(0,0)[lb]{\smash{\SetFigFont{6}{7.2}{\familydefault}{\mddefault}{\updefault}{\color[rgb]{0,0,0}alt,$J_{2n+1}$}%
}}}
\put(4691,-15766){\makebox(0,0)[lb]{\smash{\SetFigFont{8}{9.6}{\familydefault}{\mddefault}{\updefault}{\color[rgb]{0,0,0}permutation}%
}}}
\put(13951,-15176){\makebox(0,0)[lb]{\smash{\SetFigFont{7}{8.4}{\familydefault}{\mddefault}{\updefault}{\color[rgb]{0,0,0}permutation}%
}}}
\end{picture}
\end{equation} Hence we have \eqref{e10} for
  $n\ge 3$.

  Now we prove
  \begin{equation*}
    \label{e11}
    h_n(\modA ^c_{n-1}(\emptyset))\subset \BrJ.
  \end{equation*}
  It suffices to prove that every element of
  $\modA ^c_{n-1}(\emptyset)$ is a $\modZ $--linear combination of elements
  of the form depicted in the right hand side of \eqref{e12}.  
  This follows from the fact that any chord diagram on a circle is, up to $4T$ relations \cite{BN}, a linear 
  combination of chord diagrams of this form (where the dotted circle is replaced by a solid one).  
  Chmutov and Duzhin proved a quite similar statement \cite[Theorem 2.2]{CD2}, where they also assume 
  the FI relation (see \fullref{F07graph}).  
  For trivalent diagrams, a $4T$ relation can be easily derived from the IHX relation, but we do not 
  have the FI relation.  
  It is easy to avoid the latter in the proof of Chmutov--Duzhin (it is only used in one place
  \cite[Section 3.4.2.2]{CD2}) to prove the above-mentioned result on chord diagrams.  
  
  The case $n=2$ follows, since $\Br(\bJ_4(3))$ is $\modZ$--spanned by
  $[U;T_1]_{J_5}=\half[U;T_1,\tilde T_1]_{J_5}$, and
  $h_2(\modA ^c_1(\emptyset))$ is $\modZ $--spanned by $[U;T_1,\tilde
  T_1]_{J_5}$.  (Here `$1$' in `$T_1$' denotes the unit in the
  (trivial) symmetric group $S_1$ of order $1$.)

  (2)\qua Since $\xi _{k}$ is an isomorphism over $\mathbb{Q}$, the assertion follows
  immediately from (1) above and \fullref{r13}.
\end{proof}

\begin{corollary}\label{cor}
  Let $f$ be a $\modZ $--valued invariant of degree $2n$ for
  $(n{+}1)$--component links in $S^3$.  Then $f(L)$ restricted to
  Brunnian links is invariant under permutation of components and
  under orientation reversal of any component.
\end{corollary}

\bibliographystyle{gtart}
\bibliography{link}

\begin{thebibliography}{}
\providecommand\bibmarginpar{\leavevmode\marginpar}
\def\urlstyle#1{{\tt #1}}

\bibitem{BN}
\textbf{D Bar-Natan}, \href{http://dx.doi.org/10.1016/0040-9383(95)93237-2}
  {\emph{On the {V}assiliev knot invariants}}, Topology 34 (1995) 423--472
  \xox{MR}{1318886}

\bibitem{BN2}
\textbf{D Bar-Natan}, \emph{Vassiliev homotopy string link invariants}, J. Knot
  Theory Ramifications 4 (1995) 13--32 \xox{MR}{1321289}

\bibitem{BNGRT2}
\textbf{D Bar-Natan}, \textbf{S Garoufalidis}, \textbf{L Rozansky},
  \textbf{D\,P Thurston}, \emph{Wheels, wheeling, and the {K}ontsevich integral
  of the unknot}, Israel J. Math. 119 (2000) 217--237 \xox{MR}{1802655}

\bibitem{BNGRT}
\textbf{D Bar-Natan}, \textbf{S Garoufalidis}, \textbf{L Rozansky},
  \textbf{D\,P Thurston}, \href{http://dx.doi.org/10.1007/s00029-002-8109--z}
  {\emph{The {\AA}rhus integral of rational homology 3--spheres II:
  {I}nvariance and universality}}, Selecta Math. $($N.S.$)$ 8 (2002) 341--371
  \xox{MR}{1931168}

\bibitem{BLT}
\textbf{D Bar-Natan}, \textbf{T\,T\,Q Le}, \textbf{D\,P Thurston}, \emph{Two
  applications of elementary knot theory to {L}ie algebras and {V}assiliev
  invariants}, Geom. Topol. 7 (2003) 1--31 \xox{MR}{1988280}

\bibitem{CD2}
\textbf{S\,V Chmutov}, \textbf{S\,V Duzhin},
  \href{http://dx.doi.org/10.1142/S0218216594000101} {\emph{An upper bound for
  the number of {V}assiliev knot invariants}}, J. Knot Theory Ramifications 3
  (1994) 141--151 \xox{MR}{1279917}

\bibitem{CD}
\textbf{S Chmutov}, \textbf{S Duzhin},
  \href{http://dx.doi.org/10.1023/A:1010773818312} {\emph{The {K}ontsevich
  integral}}, Acta Appl. Math. 66 (2001) 155--190 \xox{MR}{1837618}

\bibitem{Cochran}
\textbf{T\,D Cochran}, \emph{Concordance invariance of coefficients of
  {C}onway's link polynomial}, Invent. Math. 82 (1985) 527--541
  \xox{MR}{811549}

\bibitem{GO}
\textbf{S Garoufalidis}, \textbf{T Ohtsuki}, \emph{On finite type 3--manifold
  invariants III: {M}anifold weight systems}, Topology 37 (1998) 227--243
  \xox{MR}{1489202}

\bibitem{Gusarov:91}
\textbf{M\,N Goussarov}, \emph{A new form of the {C}onway--{J}ones polynomial
  of oriented links}, Zap. Nauchn. Sem. Leningrad. Otdel. Mat. Inst. Steklov.
  (LOMI) 193 (1991) 4--9, 161 \xox{MR}{1157140}

\bibitem{Gusarov:94}
\textbf{M Goussarov}, \emph{On $n$--equivalence of knots and invariants of
  finite degree}, from: ``Topology of manifolds and varieties'', Adv. Soviet
  Math. 18, Amer. Math. Soc., Providence, RI (1994)  173--192 \xox{MR}{1296895}

\bibitem{Goussarov01}
\textbf{M\,N Goussarov}, \emph{Variations of knotted graphs. {T}he geometric
  technique of $n$--equiva\-lence}, Algebra i Analiz 12 (2000) 79--125
  \xox{MR}{1793618}

\bibitem{HM}
\textbf{N Habegger}, \textbf{G Masbaum},
  \href{http://dx.doi.org/10.1016/S0040-9383(99)00041-5} {\emph{The
  {K}ontsevich integral and {M}ilnor's invariants}}, Topology 39 (2000)
  1253--1289 \xox{MR}{1783857}

\bibitem{Hb}
\textbf{K Habiro}, \emph{Brunnian links, claspers, and Goussarov--Vassiliev
  finite type invariants}, to appear in Math. Proc. Camb. Phil. Soc.

\bibitem{H}
\textbf{K Habiro}, \href{http://dx.doi.org/10.2140/gt.2000.4.1} {\emph{Claspers
  and finite type invariants of links}}, Geom. Topol. 4 (2000) 1--83
  \xox{MR}{1735632}

\bibitem{HM1}
\textbf{K Habiro}, \textbf{J-B Meilhan}, \emph{Finite type invariants and
  Milnor invariants for Brunnian links} \xox{arXiv}{math.GT/0510534}

\bibitem{Kontsevich}
\textbf{M Kontsevich}, \emph{Vassiliev's knot invariants}, from: ``I. M.
  Gel'fand Seminar'', Adv. Soviet Math. 16, Amer. Math. Soc., Providence, RI
  (1993)  137--150 \xox{MR}{1237836}

\bibitem{Le}
\textbf{T\,T\,Q Le}, \emph{An invariant of integral homology 3--spheres which
  is universal for all finite type invariants}, from: ``Solitons, geometry, and
  topology: on the crossroad'', Amer. Math. Soc. Transl. Ser. 2 179, Amer.
  Math. Soc., Providence, RI (1997)  75--100 \xox{MR}{1437158}

\bibitem{LM}
\textbf{T\,T\,Q Le}, \textbf{J Murakami},
  \href{http://dx.doi.org/10.1016/S0022-4049(96)00054-0} {\emph{Parallel
  version of the universal {V}assiliev--{K}ontsevich invariant}}, J. Pure Appl.
  Algebra 121 (1997) 271--291 \xox{MR}{1477611}

\bibitem{LMO}
\textbf{T\,T\,Q Le}, \textbf{J Murakami}, \textbf{T Ohtsuki},
  \href{http://dx.doi.org/10.1016/S0040-9383(97)00035-9} {\emph{On a universal
  perturbative invariant of 3--manifolds}}, Topology 37 (1998) 539--574
  \xox{MR}{1604883}

\bibitem{lescop}
\textbf{C Lescop}, \emph{Introduction to the Kontsevich integral of framed
  tangles}, lecture notes (1999)

\bibitem{Lin}
\textbf{X-S Lin}, \emph{Power series expansions and invariants of links}, from:
  ``Geometric topology (Athens, GA, 1993)'', AMS/IP Stud. Adv. Math. 2, Amer.
  Math. Soc., Providence, RI (1997)  184--202 \xox{MR}{1470727}

\bibitem{jb}
\textbf{J-B Meilhan}, \emph{On surgery along Brunnian links in 3--manifolds}
  \xox{arXiv}{math.GT/0603421}

\bibitem{M}
\textbf{J Milnor},
  \href{http://links.jstor.org/sici?sici=0003-486X(195403)2:59:2%3C177:LG%3E2.%
0.CO%3B2--P} {\emph{Link groups}}, Ann. of Math. $(2)$ 59 (1954) 177--195
  \xox{MR}{0071020}

\bibitem{O}
\textbf{T Ohtsuki}, \href{http://dx.doi.org/10.1142/S0218216596000084}
  {\emph{Finite type invariants of integral homology 3--spheres}}, J. Knot
  Theory Ramifications 5 (1996) 101--115 \xox{MR}{1373813}

\bibitem{O2}
\textbf{T Ohtsuki}, \emph{Quantum invariants}, Series on Knots and Everything
  29, World Scientific Publishing Co., River Edge, NJ (2002) \xox{MR}{1881401}

\bibitem{thurston}
\textbf{D\,P Thurston}, \emph{Wheeling: A diagrammatic analogue of the Duflo
  isomorphism}, PhD thesis, University of California Berkeley (2000)
  \xox{arXiv}{math.QA/0006083}

\bibitem{Vassiliev}
\textbf{V\,A Vassiliev}, \emph{Cohomology of knot spaces}, from: ``Theory of
  singularities and its applications'', Adv. Soviet Math. 1, Amer. Math. Soc.,
  Providence, RI (1990)  23--69 \xox{MR}{1089670}

\end{thebibliography}

\end{document}